\documentclass{amsproc}
\usepackage{amssymb}
\usepackage{amsmath}          
\usepackage{amscd}

\newtheorem{theorem}{Theorem}[section]
\newtheorem{lemma}[theorem]{Lemma}

\theoremstyle{definition}
\newtheorem{definition}[theorem]{Definition}
\newtheorem{example}[theorem]{Example}

\newtheorem{prop}[theorem]{Proposition}
\newtheorem{coro}[theorem]{Corollary}

\newtheorem{remark}[theorem]{Remark}

\numberwithin{equation}{section}

\newcommand{\ts}{\hspace{0.5pt}}

\newcommand{\CC}{\mathbb{C}\ts}
\newcommand{\RR}{\mathbb{R}\ts}

\newcommand{\ZZ}{\mathbb{Z}}
\newcommand{\NN}{\mathbb{N}}

\newcommand{\TT}{\mathbb{T}}

\newcommand{\vL}{\varLambda}
\newcommand{\vG}{\varGamma}
\newcommand{\vOL}{\varOmega(\vL)}

\newcommand{\Lperp}{\widetilde{L}^\perp}
\newcommand{\Lnull}{L^\circ}

\newcommand{\oplam}{\mbox{\Large $\curlywedge$ }   }

\newcommand{\Oo}{(\varOmega(\vL),\alpha)}

\newcommand{\cP}{\mathcal{P}}

\newcommand{\cA}{\mathcal{A}}

\newcommand{\cCp}{\mathcal{C}_p(\vOL)}

\newcommand{\cF}{\mathcal{F}}

\newcommand{\cS}{\mathcal{S}}

\newcommand{\gammahat}{\widehat{\gamma_\vL}}

\newcommand{\Hm}[1]{\leavevmode{\marginpar{\tiny%
$\hbox to 0mm{\hspace*{-0.5mm}$\leftarrow$\hss}%
\vcenter{\vrule depth 0.1mm height 0.1mm width \the\marginparwidth}%
\hbox to
0mm{\hss$\rightarrow$\hspace*{-0.5mm}}$\\\relax\raggedright #1}}}

\newcommand{\Uphi}{U^{(\varphi) } }
\newcommand{\Pphi}{P^{(\varphi)}} 
 
\newcommand{\ani}{A^{(\varphi)}_n}
\newcommand{\aNi}{ A^{(\varphi)}_N}

\begin{document}

\title[Aperiodic order via dynamical systems]
{Aperiodic order via dynamical systems: Diffraction for sets of finite local
  complexity}

\author{Daniel Lenz}
\address{ Fakult\"at f\"ur Mathematik, D- 09107 Chemnitz, Germany  }
\curraddr{Dept. of Mathematics, Rice University, P. O. Box 1892, Houston, TX 77251}
\email{dlenz@mathematik.tu-chemnitz.de }
 \urladdr{http://www.tu-chemnitz.de/mathematik/analysis/dlenz}



\begin{abstract}  
We give an introduction into  diffraction theory for aperiodic order.  We focus on  an approach via dynamical systems and  the phenomenon of pure point diffraction.  We review recent results and sketch proofs. We then present a new  uniform Wiener/Wintner type result generalizing various earlier results of this type. 
\end{abstract}

\maketitle

\section{Introduction}
A lattice is the simplest instance of a long range ordered structure in Euclidean space.  Aperiodic order is concerned with long range ordered structures beyond lattices. A most prominent example is the Penrose tiling of the plane. While examples exhibit specific order features there is no axiomatic framework for aperiodic order yet. 

Aperiodic order has attracted a lot of attention in the last twenty five years both in physics and mathematics. One reason is the actual discovery of physical substances, later called quasicrystals, exhibiting such a form of (dis)order \cite{SBGC,INF}. Another reason is the overall interest in (dis)ordered structures. In this context aperiodic order plays a distinguished role as being situated at the border between order and disorder. 
Accordingly, various aspects have been investigated. They include geometric, combinatorial, topological and operator theoretic aspects, see e.g. the monographs \cite{Sen,Jan} and the survey and proceeding collections \cite{BMbook,Mbook,Pat,Tre}.

Here, we will deal with diffraction  i.e. harmonic analysis of aperiodic order. Diffraction is a central topic as quasicrystals were discovered by their unusual diffraction patterns. These patterns display sharp peaks indicating long range order. At the same time these patterns have five fold symmetries thereby excluding a lattice structure. In fact,  on a more conceptual level harmonic analysis of aperiodic order  had been developed, quite before the discovery of quasicrystals in work of Meyer \cite{Mey}. This work is motivated by the  question which sets allow for a "Fourier type expansion".  The corresponding sets are now known as Meyer sets and play a central role in the theory. 

Our aim here is to give an introduction into diffraction theory  of aperiodic order from the point of view of dynamical systems. This point of view has proven to be rather fruitful  as it allows to phrase both combinatorial/geometric features and Fournier analytic properties in a common framework.  It also shows clear similarities to the theory of subshifts over a finite alphabet. We therefore hope that this article can  serve as a starting point for people in dynamical systems, who are interested in aperiodic order and diffraction.  

As is clear from the size of this article, we do not intend to give a comprehensive treatment of diffraction theory. We rather focus on the phenomenon of pure point diffraction and its conceptual  understanding via dynamical systems.  In particular, we neither  discuss mixed spectra nor primitive substitutions.  

Most results covered in this article are known. We have tried to sketch proofs in a pedagogical way.  The  article also contains some new material. This concerns an observation on  symmetry in  Section \ref{Symmetry}, which seems not to be contained  explicitely in the literature. Furthermore, the uniform Wiener/Wintner type result given in Section \ref{WW} is new. It generalizes earlier results of Robinson \cite{Rob}, Walters \cite{Wal2} and Lenz \cite{Le}.  Our proof  follows the method given in  \cite{Le}. 

The paper is organized as follows: In Section \ref{Point} we introduce the point sets  of interest and the associated dynamical systems. Section \ref{Diffraction} gives an introduction into diffraction theory. The main results as well as their history are discussed in Section \ref{Results}. Proofs are sketched in Section \ref{Proofs}. Section \ref{Symmetry} contains the observation on how symmetries of points sets show up in the corresponding diffraction.  Meyer sets  and more specially regular model sets are studied  Section \ref{Meyer}. As Section \ref{Proofs} shows, uniform Wiener/Wintner type results are useful in the study of diffraction. Thus, we present our  new result of this form in Section \ref{WW}.
Finally, Section \ref{Further} contains further remarks and open questions.

\section{Point sets with finite local complexity and the associated dynamical systems}\label{Point}
Point sets with finite local complexity can be seen as geometric analogues of
sequences taking only finitely many values.  The associated dynamical systems
are geometric analogues of subshifts over a finite alphabet.  This point of view has been developed over the last fifteen years or so. In this section, we give an introduction into this topic.

\bigskip

Our basic setup is as follows: We consider subsets of Euclidean space $\RR^N$. The Euclidean norm is denoted
by $\| \cdot\|$ and the closed ball around the origin $0$ with radius $S$ is denoted by
$B_S$.  The Lebesgue measure of a measurable subset of $\RR^N$ is denoted by $|M|$ and the cardinality of a set $F$ is denoted by $\sharp F$.

\begin{definition} Let $\vL$ be a subset of $\RR^N$. Then,  $\vL$ is called
  uniformly discrete if there exists $r>0$ with
$$\|x-y\|\geq 2 r$$
 for all
  $x,y\in \vL$ with $x\neq y$.  The set $\vL$ is called relatively
  dense if  there exists an $R>0$ with
$$\RR^N =\cup_{x\in \vL} ( x + B_R ).$$
  If $\vL$ is both uniformly discrete (with parameter $r$)  and
  relatively dense (with parameter $R$)  it is
  called a Delone set or an  $(r,R)$ - Delone set. 
\end{definition}

\begin{remark} If  $\vL$ is uniformly discrete with parameter $r$, then  open
balls around  points of $\vL$ with radius $r$ are disjoint.  This is the  reason for the factor $2$ appearing in the above definition.  The largest $r$
with this property is called the packing radius of $\vL$.  On the other hand
if  $\vL$ is relatively dense with parameter $R$, then no point of $\RR^N$ has distance larger than $R$ to $\vL$. Then the smallest $R$  with with property is  called the covering radius of $\vL$. 
\end{remark}

\medskip

We now introduce the crucial concept of patch. A patch is a local
configuration in a Delone set. Various versions are considered in the
literature. For our purposes the following seems the most practical. 
A patch   of size $S>0$ in a Delone set $\vL$ is a set of the form
$$ (\vL - x)\cap B_S,$$
where $x$ belongs  to $\vL$. Thus, any  patch contains the origin.  Sometimes these patches are called centered ball patches. We define
$$N_\vL (S):= \sharp \{ (\vL -x)\cap B_S : x\in \vL\}.$$

We are interested in Delone sets whose patches satisfy a certain finiteness
condition.  This condition is characterized next.

\begin{lemma} Let $\vL$ be a $(r,R)$ Delone set. The following assertions are
  equivalent:

\begin{itemize}
\item[(i)] For each $S>0$, the number $N_\vL (S)$ is finite, i.e. there  are only finitely many patches of size $S$ in
  $\vL$. 
\item[(ii)] The set $\vL -\vL$ is discrete and closed.
\item[(iii)] The set $(\vL - \vL )\cap B_S$ is finite for any $S>0$.
\item[(iv)] The number $N_\vL (2 R)$ is finite.   
\end{itemize}
\end{lemma}

The equivalence between $(i)$, $(ii)$ and $(iii)$ is straightforward. The equivalence
of $(i)$ and $(iv)$ is due to  Lagarias, see Corollary 2.1 of \cite{Lag}.

\begin{definition} Let $\vL$ be a Delone set. Then $\vL$ is said to have
  finite local complexity (FLC), if it satisfies one of the conditions of the
  previous lemma. 
\end{definition}

By condition (iv) of the previous lemma,  Delone sets with finite local
complexity can be considered as   geometric analogues of  sequences over a
finite alphabet.  In fact, it is easily possible to  associate one
dimensional Delone sets with (FLC) to sequences over a finite alphabet and
vice versa. This is discussed in some detail in the next example.

\medskip

\begin{example} Let $\cA$ be a finite set. To each $a\in \cA$ associate a
finite interval $[0,l_a]$ in $\RR$ by choosing $0<l_a < \infty$.  Then, we can  obtain a Delone set $\vL_\omega$ for  any sequence
$\omega:\ZZ\longrightarrow \cA$  by ``tiling'' $\RR$ with the  intervals
$[0,l_{\omega(n)}]$ in the obvious way  according to 
$$    \cdots    [0, l_{\omega(-1)}] |     [0,l_{\omega(0)}] [0,l_{\omega(1)}]\cdots,$$
where $|$ denotes the position of the origin. More precisely, 
$$\vL_\omega:=\{0\}\cup\{ \sum_{j=0}^n l_{\omega(j)} : n\geq 0\}\cup\{ - \sum_{j= -n}^{-1} l_{\omega(j)} : n \geq 1\}.$$
The Delone set $\vL_\omega$  has (FLC). It contains the origin and if the lengths $l_a$,
$a\in \cA$, are pairwise different, we can recover $\omega$ from
$\vL_\omega$.

Conversely a one
dimensional (r,R)- Delone set containing the origin  with  FLC can be
converted into a sequence with values in the finite set $\{(\vL - x)\cap B_{2
  R} : x\in \vL\}$ as follows: Enumerate the points of $\vL$  in increasing
order according to
$$  \cdots < x_{-1} <x_0=0< x_1< x_2 <\cdots$$
and then define $\omega_\vL: \ZZ\longrightarrow  \{(\vL - x)\cap B_{2
  R} : x\in \vL\}, \omega_\vL (n):= (\vL - x_n)\cap B_{2R}.$
From this sequence we can then recover $\vL$. 
\end{example}

\medskip

These considerations show that from a combinatorial point of view one
dimensional sets with (FLC) and sequences over a finite alphabet are
essentially equivalent.  This means, in particular, that  $(FLC) $ is not really an order
requirement as any sequence (no matter how disordered it is) gives rise to a
Delone set with (FLC). The example also shows that in general Delone sets with
(FLC) will not have the  property that $\vL -\vL$ is uniformly discrete. To
see this is suffices to consider $\cA=\{0,1\}$ and $l_0=1$ and $l_1=\alpha$
with $\alpha$ irrational. Then, any "typical" sequence $\omega:
\ZZ\longrightarrow \cA$ will give rise to a Delone set $\vL_\omega$ whose set
of differences is not uniformly discrete.  
 
The equivalence of
one dimensional sets with (FLC) and sequences over a finite alphabet breaks
down when  it comes to comparing the associated  dynamical systems.  This  is studied in work of Clark/Sadun \cite{CS1} (see also their work \cite{CS2} for higher dimensional analogues).

\medskip

We will now discuss  two  regularity properties that a Delone  may have. In
order to formulate them, we introduce the concept of locator set. 
The locator set  $L_\vL^P $ of the patch $P$ of size $S$ in $\vL$ is the set of all points
in $\vL$ at which $P$ occurs, i.e. 
$$ L_\vL^P :=\{ x\in \vL : (\vL -x)\cap B_S =P\}.$$

\begin{definition} A Delone set $\vL$ is said to be repetitive if $L_\vL^P$ is
  relatively dense for any patch $P$ of $\vL$.
\end{definition}

The other property can be described in various ways:

\begin{lemma} Let $\vL$ be Delone and $P$ a patch in $\vL$. Then, the following assertions are equivalent:
\begin{itemize}
\item[(i)] For any sequence $(p_n)$ in $\RR^N$ the limit
$\lim_{n\to \infty} \frac{\sharp L_\vL^P \cap (p_n + B_n) }{ |B_n|  }$
exists. 
\item[(ii)] There exists a number $\nu_P$ such that for any $\varepsilon >0$ there exists an $S>0$ with  $| \nu_P - \frac{\sharp L_\vL^P \cap (p + B_S ) }{ |B_S|}| \leq \varepsilon$ for all $p\in \RR^N$. 
\end{itemize}
\end{lemma}
\begin{proof} The implication $(ii)\Longrightarrow (i)$ is clear. As for $(i)\Longrightarrow (ii)$,  interspersing sequences shows that the limits in $(i)$ must be independent of the sequence $(p_n)$. Now, $(ii)$ follows easily. 
\end{proof}

\begin{definition} A Delone set $\vL$ is said to have uniform patch
  frequencies (UPF) if for any patch $P$ in $\vL$ one of the conditions of the previous lemma holds.  The  number $\nu_P$ is then called the frequency of $P$. 
\end{definition}

We are now heading towards introducing dynamical systems associated to Delone
sets. To a discrete set  $\vL$ let $\cP(\vL)$ be the set of all patches of $\vL$,
i.e. 
$$\cP (\vL):=\{ (\vL - x)\cap B_S : x\in \vL, S>0\}.$$
Then, we define the hull $\vOL$  of the Delone set $\vL$ by
$$ \vOL:=\{ \vG \subset \RR^N: \vG\neq \emptyset \:\mbox{and}\:  \cP (\vG)\subset \cP (\vL)\}.$$
If  $\vL$ is an $(r,R)$ - Delone so must be any $\vG\in \vOL$ by construction. 
Obviously, $\vOL$ is invariant under translations. Thus, we have an action
$\alpha$  of
$\RR^N$ on $\vOL$ by
$$ \alpha  :\RR^N \times \vOL\longrightarrow \vOL, \alpha_t (\vG):= t + \vG.$$
We will equip the set of all Delone sets with a metric.
The basic idea behind the metric is that Delone sets are close whenever they
agree on a large ball around the origin up to a small translation. To make
this precise, we set for $\vG,\vL$ Delone
$$ \widetilde{d}(\vL,\vG):=\inf\{\varepsilon >0 : \exists x,y\in B_\varepsilon\;
\mbox{s.t} \: (\vL - x) \cap B_{1/\varepsilon} = (\vG- y)\cap
B_{1/\varepsilon}\}.$$
Note that the infimum is finite, as the two sets in question are not empty. We
then define  for  $\vG,\vL$ Delone
$$ d(\vG,\vL):=\min\{\frac{1}{\sqrt{2}}, \widetilde{d}(\vG,\vL)\}.$$
Then, $d$ is a metric. Symmetry and non-degeneracy are clear. The cut-off with
$1/\sqrt{2}$ makes $d$
satisfies the triangle inequality as well (see e.g.  \cite{LMS}). This metric defines a topology. 
Convergence with respect to this topology can directly be seen to have the following properties. 

\begin{lemma} Let $\vG$ be a $(r,R)$ - Delone set containing the origin and $(\vG_n)$ be a sequence of Delone sets converging to $\vG$. 

(a) There exists a sequence $(t_n)$ in $\RR^N$ converging to $0$ such that  $(\vG_n -t_n)$ contains the origin for each $n$ and $(\vG_n - t_n)$ converges to $\vG$.  

(b) If each $\vG_n$ contains the origin, then there exists for any $S>0$ an $N$ with $\vG_n \cap B_S = \vG \cap B_S$ for all $n\geq N$. 
\end{lemma}

It is not hard to see that $\alpha$ defined above is a continuous action on
the set of all Delone sets. Hence, each $\vL$ gives rise to a topological
dynamical system $\Oo$. 
 There is  a ``dictionary'' between basic  properties of $\vL$ and basic
properties of $\Oo$ given in the next three theorems. While these results are well known we include sketches of proofs for the convenience of the reader. Here, we are concerned with Delone sets in Euclidean space. The results below also hold for Delone sets on locally compact Abelian groups. This is investigated by Schlottmann in \cite{Schl}. 

We start with equivalence of (FLC) and compactness of $\vOL$.  The corresponding result for symbolic dynamics is a direct consequence of Tychonoffs Theorem. 
For tilings the result is due to Radin/Wolff  \cite{RW}.  For Delone sets it  can be found in \cite{LP}.

\begin{theorem} Let $\vL$ be a Delone set. Then, $\vL$ has (FLC) if and only
  if its hull $\vOL$ is compact. 
\end{theorem}
\begin{proof}
  Let $\vL$ be an $(r,R)$ - Delone set.  To show that (FLC) implies
  compactness, consider a sequence $(\vG_n)$ in $\vOL$. We have to provide a
  converging subsequence.  Each element of $\vOL$
  contains a point in $B_R$. Hence, each $\vG_n$ contains a point in $B_R$.  These points must have an accumulation point. Shifting if necessary, we can assume 
  without loss of generality that this accumulation point is $0$. By
  shifting and going to a subsequence we can then assume without loss of
  generality that each $\vG_n$ contains the origin. For any $k\in \NN$, we 
  then consider the set $\{\vG_n \cap B_k: n\in \NN\}$. By (FLC) this set is
  finite.  By doing a diagonal sequence argument we conclude the desired
  statement.

  Conversely, let $\vOL$ be compact.  If $\{(\vL -x)\cap B_S : x\in \vL\}$
  were infinite for some $S>0$, we could find $x_j\in \vL$, $j=1, 2,  \ldots$ such
  that $(\vL -x_j)$ are pairwise different on $B_S$ and, obviously, all contain the origin.  Then, the sequence $(\vL
  - x_j)$ can not have an accumulation point.
\end{proof}

We now come to a characterization of  repetitivity. For symbolic dynamics the result is well known (see e.g. \cite{Que}).  For tilings it can be found in \cite{Sol}. For Delone sets it is discussed in \cite{LP}.

\begin{theorem} Let $\vL$ be a Delone set with (FLC). Then, $\vL$ is
  repetitive if and only if $\Oo$ is minimal (i.e. each orbit is dense). 
\end{theorem}
\begin{proof} 
  Let $\Oo$ be minimal. If $\vL$ were not repetitive, we could find
  arbitrarily large balls on which a certain patch $P$ does not occur. Shifting
  these balls to the origin and using compactness, we would obtain a $\vG$
  which would not contain $P$ at all.  Therefore, translates of $\vG$ could not approximate $\vL$.  This contradicts minimality.
  
  Conversely, let $\vL$ be repetitive.  Let $P$ be an arbitrary patch in $\vL$.  Then,
  there exists an $S>0$ such that any ball of size $S$ in $\vL$ contains a
  translate of $P$.  Hence, any $\vG\in \vOL$ must contain a copy of $P$ in
  the ball of size $S$ around the origin.  As $P$ is arbitrary, minimality
  follows.
\end{proof}

We finally discuss equivalence of (UPF) and unique ergodicity. Again, this is
well known for symbolic dynamics (see e.g. the books \cite{Que,Wal}). For tilings it is discussed in \cite{Sol} and for Delone sets  in \cite{LP,LMS}. 

\begin{theorem} Let $\vL$ be a Delone set with (FLC). Then, $\vL$ has  uniform
  patch frequencies (UPF) if and only if $\Oo$ is uniquely ergodic  (i.e.
  there exists a unique invariant probability measure on $\vOL$). 
\end{theorem}
\begin{proof} 
  It is well known that unique ergodicity is equivalent to uniform existence
  of the limits in Birkhoff ergodic theorem for a sufficiently large set of
  continuous functions. It turns out that patterns can be used to create such
  a set. More precisely, define for a pattern $P$ and $\varphi \in C_c
  (\RR^N)$ the function $f_{P,\varphi}$ on $\vOL$ by
$$ f_{P,\varphi} (\vG) = \sum_{x\in L_\vG^P} \varphi (-x).$$
Note that the sum has only finitely many non vanishing terms as $\varphi$ has
compact support. 
The sign in $-x$ does not play role and is only to make this
consistent with later considerations. These functions are continuous
functions. UPF can be seen to be equivalent to uniform convergence of  the Birkhoff averages for these functions. This then turns out to be equivalent to unique
ergodicity. 
\end{proof}

Let us finish this section by recalling some basic facts on spectral theory of
dynamical systems.  Let $\vL$ be Delone  with (FLC) and $m$ an $\alpha$-invariant measure
on $\vOL$. The action  $\alpha$ on $\vOL$  then induces a unitary representation $T$
of $\RR^N$ on $L^2 (\vOL,m)$ viz
$$ (T_t f) (\vG) = f(- t + \vG).$$
An $f\in L^2 (\vOL,m)$ with $f\neq 0$  is called an eigenfunction (to the
eigenvalue $\xi$) if 
$$ T_t f = \exp (i \xi t) f$$
for all $t\in \RR^N$ (where the equality is understood in  the $L^2$
sense). An eigenvalue is called a continuous eigenvalue if it
admits a continuous eigenfunction $f$ satisfying
$$ f(-t + \vG) = \exp(i \xi t) f(\vG)$$
for all $t\in \RR^N$ and all $\vG\in \vOL$.  $\Oo$ is said to have pure point
spectrum if $L^2 (\vOL,m)$ has a basis consisting of eigenfunctions.

\section{Diffraction theory}\label{Diffraction}
In this section we present a basic setup for diffraction \cite{Cow}.  For models  with aperiodic order this  framework has  been advocated by Hof \cite{Hof}  and become a standard by now. The crucial quantity is a measure, called the diffraction measure and denoted by  $\gammahat$. This measure  represents the intensity (per unit volume). It models the outcome  of a diffraction experiment.

\bigskip

In a diffraction experiment a solid  is put into an  incoming beam of e.g. $X$ rays. The atoms of the solid then interact with the beam and one obtains an outcoming wave. The intensity of this wave  is then measured on a screen.  When modeling diffraction, the two basic principles are the following:

\begin{itemize}
\item Each point $x$ in the solid gives rise to a wave $\xi \mapsto \exp(-i x \xi)$. The overall wave $w$ is the sum of the single waves. 

\item The quantity measured in an experiment is the intensity given as the square of the modulus of the wave function. 
\end{itemize}

We start with  by implementing this for a finite set $F\subset \RR^N$.  Each $x\in F$ gives rise to a wave
$\xi \mapsto \exp(-i x \xi)$. The overall wavefunction $w_F$  induced by $F$ is then
$$ w_F (\xi)  = \sum_{x\in F} \exp(-i x \xi).$$
Thus, the intensity $I_F$ is the function  given as
$$I_F (\xi) = \sum_{x,y\in F} \exp(-i (x-y) \xi) = \cF (\sum_{x,y\in F}
\delta_{x-y}).$$
Here, $\delta_z$ is the unit point mass at $z$ and $\cF$ denotes the Fourier
transform.  To describe diffraction for a solid with many atoms it is common
to model the solid by a Delone set in $\RR^N$.  When trying to establish a formalism  as above for an infinite set $\vL$,
one faces the immediate problem that
$$ w_\vL =  \sum_{x\in \vL} \exp(-i x\xi)$$
does not make sense. One may try and give it a sense as a tempered
distribution. This, however,  does not solve the problem as the quantity we
are after is the intensity given as $|w_\vL|^2$. Now, neither modulus nor
products are defined for distributions. This is not only a mathematical
issue. There is a physical reason  behind the divergence: The intensity of the
whole set is infinite. The correct quantity to consider is not the intensity
but a normalized intensity viz.  the intensity per unit
volume. We therefore try and define
$$
I =\lim_{n\to\infty} \frac{1}{|B_n|} I_{\vL\cap B_n}.$$
Various comments
are in order: As $\vL$ is uniformly discrete, $\vL\cap B_n$ is finite.
Thus, $I_{\vL\cap B_n}$ is defined. Thus, on the right hand side we have a
sequence of functions. We consider this sequence as a sequence of measures by
considering each function as the density with respect to Lebesgue measure. The
limit is then meant in the sense of vague convergence of measures. Recall that a sequence $(\nu_n)$  of measures converges in the vague topology
to the measure $\nu$  if $\nu_n (\varphi)\longrightarrow \nu(\varphi)$, $n\to
\infty$,  for each  continuous function  $\varphi :\RR^N\longrightarrow
\RR$ with compact support.  Of course,
it is not clear (and will be wrong in general) that the limit exists. If it
exists it is a measure. Let us emphasize  once more that this measure is the crucial
object as it describes the outcome of a physical diffraction experiment.

In order to discuss issues related to existence of the limit a little closer, we need some preparation.  Besides the concept of vague convergence of measures, which we have just defined, we will need the  Schwarz space $\cS$.  This is the space of all  functions $\varphi :\RR^N\longrightarrow \RR$,  which are infinitely many often differentiable and all of whose derivatives of any order  go faster to zero than any polynomial tends to infinity.  Moreover,  we also recall definitions concerning convolutions. 
For $\varphi, \psi\in C_c (\RR^N)$ 
we define the convolution  $\varphi \ast \psi \in C_c (\RR^N)$ by
$$\varphi \ast \psi (x) = \int_{\RR^N} \varphi (x -y) \psi (y) d y $$
and $\widetilde{\varphi}\in C_c (\RR^N)$ by $\widetilde{\varphi} (x) =
\overline{\varphi}(-x)$. The convolution of $\varphi \in C_c (\RR^N)$ with  a
measure $\nu$ on $\RR^N$ is the continuous function defined by 
$$\nu \ast \varphi (t) = \int_{\RR^N} \varphi (t -s) d\nu (s).$$
Finally, for a function $\varphi \in C_c (\RR^N)$ we define $\widetilde{\varphi} (x)= \overline{ \varphi (-x) }$.

\begin{prop} Let $\vL$ be a Delone set. The following assertions are
  equivalent:
\begin{itemize}
\item[(i)] The measures $I_{\vL\cap B_n}$,$ n\in \NN$, converge to a limit $I_\vL$ in
  the vague topology.
\item[(ii)]  The measures $\gamma_\vL^n:= \frac{1}{|B_n|} \sum_{x,y\in
    \vL\cap B_n} \delta_{x-y}$, $n\in \NN$, converge to a limit, $\gamma_\vL$, in the vague topology.
\end{itemize}
In this case, $I_\vL$ is a positive measure and the Fourier transform of
$\gamma_\vL$  in
the sense that $I_\vL (|\cF(\varphi)|^2) = \gamma_\vL (\varphi\ast \widetilde{\varphi})$ for any
$\varphi$ in  $C_c (\RR^N)$. 
\end{prop}
\begin{proof}

 A direct calculation shows that  the measure $I_{\vL\cap B_n}$ is the Fourier transform of
  $\gamma_\vL^n$ in the sense of tempered distributions i.e. 
$$ \int I_{\vL\cap B_n}(\xi) \varphi (\xi) d\xi = \gamma_\vL^n
(\cF^{-1} (\varphi))$$ 
for any function $\varphi$ in the Schwarz space.   This shows the desired equivalence in the sense of convergence of tempered distributions. 
Now, for the measures in
question convergence with respect to the vague topology is equivalent to
convergence as tempered distributions. To show this  requires some care. We first note that the measures $\gamma_\vL^n$ are uniformly translation bounded (i.e. there exists a $C$ with $\gamma_\vL^n ( t + B_1) \leq C$ for any $t\in \RR^N$ and $n\in \NN$).  Therefore, we  can replace $C_c (\RR^N)$ by $\mathcal{S}$ as far as  $\gamma_\vL^n$ are concerned. Now, note that
$\gamma_\vL^n$ are also  positive definite (i.e. $\gamma_\vL^n \ast \varphi \ast \widetilde{\varphi}(0) \geq 0$ for all $\varphi \in C_c (\RR^N)$).  This gives that the measures 
$I_{\vL\cap B_n}$ are then uniformly translation bounded as well.  Hence, we can replace $C_c (\RR^N)$ by $\mathcal{S}$ when dealing with $I_{\vL\cap B_n}$. 
  This show the equivalence between
$(i)$ and $(ii)$. 

The last statement is obvious for $\varphi \in \mathcal{S}$ and follows for $\varphi \in C_c (\RR^N)$ by approximation.  
\end{proof}

These considerations lead to the following definition.

\begin{definition} Let $\vL$ be a Delone set. The set  $\vL$ is said to have a
  well defined autocorrelation if $\frac{1}{|B_n|} \sum_{x,y\in
    \vL\cap B_n} \delta_{x-y}$ converge. The limit $\gamma_\vL$ is called the
    autocorrelation function (even though it is a measure). 
 In this case, the  Fourier transform  $\gammahat$ of $\gamma_\vL$  is called the
  diffraction measure. 
\end{definition}

\medskip

\begin{remark}  To understand the averaging in the definition of $\vL$, it is instructive to proceed as follows: Define the Dirac comb $\delta_\vG$ of the set $\vG\subset \RR^N$ by $\delta_\vG = \sum_{x\in \vG} \delta_x$. Then, 
$$\gamma_\vL = \lim_{n\to \infty} \frac{1}{|B_n|}  \delta_{ \vL \cap B_n} \ast \delta_{ - \vL \cap B_n} =\lim_{n\to \infty} \frac{1}{|B_n|}  \delta_{ \vL \cap B_n} \ast \delta_{ - \vL }.$$
Here, the approximants are convolutions of the (uniformly  in $n$) bounded measures $\frac{1}{|B_n|}  \delta_{ \vL \cap B_n} $ and the (uniformly in $n$)  translation invariant measures $\delta_{ \vL \cap B_n}$. In particular, $\gamma_\vL$ is a translation bounded measure of infinite total mass. 
\end{remark}

\bigskip

We are particularly interested in the point part of $\gammahat$. We introduce the following notation. The points $\xi\in \RR^N$ with $
\gammahat(\{\xi\})\neq 0$ are called Bragg peaks. The value
    $\gammahat( \{\xi\} )$ is called the intensity of the Bragg peak.

Let us now shortly summarize our approach so far. 
We have presented an abstract framework to deal with an diffraction
experiment. The outcome of a diffraction experiment is described by a measure,
the so called diffraction measure, $\gammahat$. In this context the following question
arise naturally:

\begin{itemize}
\item When does $\gammahat$ exist?
\item When is $\gammahat$ a pure point measure?
\item Where are the Bragg peaks?
\item Which are the intensities of the Bragg peaks?
\end{itemize}
In the next section we will   present answers to these questions  in the framework
of dynamical systems. 

\section{Results on diffraction}\label{Results}
In this section we present some answers to  the questions raised at the end of
the last section.  These answers are formulated in terms of the dynamical
system associated to $\vL$. In this sense, they can be considered as an
extension to diffraction  of the ``dictionary'' between properties of $\vL$ and properties of
$\vOL$. This section is devoted to  statements of  results
and a discussion of the literature. The next section provides some ideas for 
the proofs.

\bigskip

Throughout we consider    $\vL$ Delone satisfying  (FLC) and (UPF). 
By the considerations above this  implies in particular that   $\vL-\vL$ is
  discrete and closed,  $\vOL$ is
  compact and $\Oo$ is uniquely ergodic.  Let $m$ be the unique translation
  invariant probability measure on $\Oo$.

\medskip

The first result answers the question of existence of $\gamma_\vL$. 
\begin{theorem}\label{closedformula} Let $\vL$ be a Delone set with (FLC) and (UPF). Then,
  $\gamma_\vG$ exists for every $\vG$ in $\vOL$ and equals $\gamma_\vL$. The
  measure $\gamma_\vL$ is supported on $\vL -\vL$ and given by the closed formulas 
$$\gamma_\vL  (\varphi) = \sum_{z\in \vL - \vL}  c_z \, \varphi (z) = \int_{\vOL}
\sum_{x,y\in \vG} \sigma (x)  \varphi (x-y) d m (\vG), $$
where $c_z:= \lim_{n\to \infty} \frac{1}{|B_n|} \sharp\{x\in \vL \cap B_n : x+z \in \vL\}$ and  $\sigma\in C_c (\RR^N)$ is arbitrary with $\int \sigma (t) dt =1$.
\end{theorem}

Existence of $\gamma_\vG$ and the first close formula goes back to Hof \cite{Hof}. It has then be extended to various other contexts and situations. In particular, $\RR^N$ can be replaced by a $\sigma$-compact locally compact Abelian group \cite{Schl}. Moreover, neither unique ergodicity nor FLC are needed to obtain a closed formula for $\gamma_\vL$. In fact, it is possible to give a closed formula in the context of point processes and Palm measures \cite{Gou1} or in the context of translation bounded measures on $\sigma$-compact locally compact Abelian group \cite{BL}. The second closed formula given above is taken from \cite{BL}. 

\bigskip

We now discuss an answer to the question whether $\gammahat$ is a pure point measure. 

\begin{theorem} \label{characterization} Let $\vL$ be a Delone set with (FLC) and (UPF). The following
  assertions are equivalent:
\begin{itemize}
\item[(i)] $\gammahat$ is a pure point measure.
\item[(ii)] The dynamical system $\Oo$ has pure point dynamical spectrum. 
\end{itemize}
In this case the group of eigenvalues is  the smallest subgroup of $\RR^N$ containing all  $\xi$ with $ \gammahat\{\xi\}\neq 0$.
\end{theorem}
For symbolic dynamics this type of result has been proven by Qu\'{e}ff\'{e}lec in \cite{Que}. For Delone dynamical systems the implication $(ii)\Longrightarrow (i)$ has been shown by Dworkin \cite{Dwo} and the corresponding reasoning is known as Dworkin argument. The equivalence given above is due to Lee/Moody/Solomyak \cite{LMS}.  Their result can be extended to  rather general point processes in $\RR^N$ using Palm measures as shown by Gou\'{e}r\'{e} \cite{Gou2}. Their result can also be extended to translation bounded measures on locally compact Abelian groups, as shown by Baake/Lenz \cite{BL}.  The statement on the eigenvalues is implicit in \cite{LMS}. It can be found explicitely in \cite{BL}. 

 The argument of  \cite{Dwo} shows essentially  that the diffraction spectrum is part of the dynamical spectrum.  At the same time  there is also work of van Enter/Mi\c{e}kisz \cite{EM} showing that the dynamical spectrum is in general strictly larger than the diffraction spectrum. More precisely, they give an example  of a system whose dynamical spectrum contains both a point component and a continuous component but  the diffraction measure has only a continuous component.

\bigskip

The previous theorem can be used to show the following. Recall that $N_\vL (S)$ is the
number of patches of size $S$ in $\vL$. 

\begin{theorem}  \label{entropy} Let $\vL$ be a Delone set with (FLC) and (UPF) and  $\gammahat$
  a pure point measure. Then, the patch counting entropy of $\vL$ vanishes
  i.e. 
$$ 0=\lim_{S\to \infty} \frac{ \ln N_\vL (S)}{|B_S|}. $$
\end{theorem}

This result is due to Baake/Lenz/Richard \cite{BLR}. It confirms the intuition that  long range order (as expressed by pure point diffraction) implies  order in terms of bounds on the growth of complexity.

\bigskip

We finally come to intensity of the Bragg peaks.  The basic idea is that the intensities of Bragg peaks can be calculated via averaged Fourier coefficients.  More precisely, define for a Delone set $\vG$,  $\xi\in \RR^N$ and $S>0$ 
$$ c_S^\xi (\vG):=\frac{1}{|B_S|}  \sum_{x\in \vG\cap B_S} \exp (-i\xi x).$$
Since the work of Bombieri/Taylor a basic assumption has been that 
$$\gammahat(\{\xi\})=\lim_{S\to \infty}  |c_S^\xi (\vG)|^2.$$
This assumption has then be called Bombieri/Taylor conjecture. It was shown to hold for regular models sets by Hof \cite{Hof} and in a more general context by Schlottmann \cite{Schl} and for primitive substitutions by G\"ahler/Klitzing in \cite{GK}. The work of Hof hints at a connection to continuity of eigenfunctions.  This has been confirmed recently by Lenz \cite{Le}. There one can find the following result.

\begin{theorem}\label{bombieritaylor}
Let $\vL$ be a Delone set with (FLC) and (UPF). Assume that
  $\gammahat$ is a pure point measure and all Bragg peaks are continuous
  eigenvalues. Then, 
$$ \gammahat(\{\xi\})=\lim_{S\to \infty}  |c_S^\xi (\vG)|^2$$
for all $\vG\in \vOL$ and all $\xi\in \RR^N$.  
\end{theorem}

This results  allows one to recover the mentioned results of Hof and G\"ahler/Klitzing.  In fact, the considerations in \cite{Le} treat various further examples.

\section{Ideas of the proofs}\label{Proofs}
In this section we  sketch proofs  of the results of
the previous section.  Throughout we assume that $\vL$ is Delone with (UPF) and (FLC).

\begin{proof}[Proof of Theorem \ref{closedformula}]  By (UPF) the frequency
$$\lim_{n\to \infty} \frac{1}{|B_n|} \sharp\{x\in \vG \cap B_n : x+z \in
\vG\}$$
exists for any $z\in \vL-\vL$ and $\vG$ in $\vOL$ and is independent of $\vG$. This shows
existence of $\gamma_\vG$, its independence of $\vG$, and the first equality.
It remains to show 
$$\gamma (\varphi) =  \int \sum_{x,y\in \vL} \sigma (x)  \varphi (x-y) d\mu,
(\omega)$$
for all $\varphi \in C_c (\RR^N)$. Fix $\varphi\in C_c (\RR^N)$. As $m$ is
translation invariant, the map 
$$\sigma\mapsto  \int \sum_{x,y\in \vL} \sigma (x)  \varphi (x-y) d\mu (\omega)$$
can easily be seen to provide a translation invariant measure on  $\RR^N$.  As
there is (up to multiples) only one translation invariant measure on $\RR^N$,
we infer independence  of $\sigma$ provided $\int \sigma (t) dt =1$.  In fact,
we are also  allowed to chose functions of the form $\frac{1}{|B_S|}
\chi_{B_S}$, where $\chi$ denotes the characteristic function. Choosing such
functions, letting $S\to \infty$ and applying the ergodic theorem, we obtain
the desired equality. 
\end{proof}

We will now discuss a connection between diffraction and
spectral theory of the associated dynamical system.  This connection can be found in the work of  Dworkin \cite{Dwo} (see \cite{EM} for strongly related ideas as well). The
measure $m$ is the unique invariant probability measure on $\vOL$. $T$ denotes
the unitary representation of $\RR^N$ on $L^2 (\vOL,m)$. The inner product on
$L^2 (\vOL,m)$ is denoted by $\langle\cdot,\cdot \rangle$.  

By the Stone/von Neumann theorem, each $f\in L^2 (\vOL,m)$ gives rise to the
spectral measure $\rho_f$ on $\RR^N$. This measure is characterized by 
validity of 
$$ \langle f, T^t f\rangle = \int \exp(i t\xi)d\rho_f (\xi)$$
for all $t\in \RR^N$. The spectral measures determine the whole spectral
theory of $T$. In particular, a spectral measure is a pure point measure if
and only if the corresponding function is a linear combination of
eigenvectors. Thus,  $T$ has pure point spectrum if and only if all
$\rho_f$ are pure point measures.

Each $\varphi \in C_c (\RR^N)$ induces a continuous function $f_\varphi$ on
$\vOL$ given by
$$ f_\varphi (\vG) = \sum_{x\in \vG} \varphi (-x).$$

The connection between diffraction spectrum and the dynamical spectrum is then
given by the following lemma.

\begin{lemma}\label{link} For $\varphi,\psi \in C_c (\RR^N)$ the equality
$$\gamma_\vL \ast \varphi\ast \widetilde{\psi} (t) =\langle f_\varphi, T^t
f_\psi \rangle$$
holds for all $t\in \RR^N$. In particular, 
$$ |\cF (\varphi)|^2 \gammahat = \rho_{f_\varphi}.$$
\end{lemma}
\begin{proof} The first statement can be derived from the second closed formula in
  Theorem \ref{closedformula}   by a direct  but somewhat lengthy computation \cite{BL}. The second
  statement then follows by taking Fourier transforms. 
\end{proof}

Having discussed  this connection we can now sketch proofs for the results of
the previous section. 

\begin{proof}[Proof of Theorem \ref{characterization}]
We start by discussing the equivalence between $(i)$ and $(ii)$. 

$(i)\Longrightarrow (ii)$: If $T$ has
pure point spectrum, then certainly all $\rho_f$, $f\in L^2 (\vOL,m)$, are pure
point measures. Hence, by Lemma \ref{link}, all measures of the form $|\cF
(\varphi)|^2 \gammahat$ are pure point measures. Hence, $\gammahat$ is a pure
point measure.

$(ii)\Longrightarrow (i)$: Let  $\gammahat$ be a pure point measure. Consider
the set $\cCp$ consisting of all continuous functions on $\vOL$ whose
spectral measure is a pure point measure.  We show that this set is an
algebra, which contains the constant functions,  is closed under complex
conjugation and separates the points:  The constant functions belong
to $\cCp$, as they are eigenvectors to the eigenvalue $0$. As the complex conjugate of an eigenfunction is an
eigenfunction,  $\cCp$ is closed under complex conjugation. As 
the product of two (bounded) eigenfunctions is an
eigenfunction, it is possible to show that $\cCp$ is closed under products \cite{LMS,BL}. 
Moreover, by Lemma \ref{link} again,  all spectral measures $\rho_{f_\varphi}$
are pure point measures. This implies that all $f_\varphi$ belong to
$\cCp$. These $f_\varphi$ obviously separate the points of $\vOL$.

These considerations show that $\cCp$ satisfies the assumptions of
Stone/Weierstrass Theorem. Hence, we conclude that $\cCp$ is dense (with
respect to the supremum norm) in the continuous function on $\vOL$. Then, it
must  also be dense (with respect to the $L^2$ norm) in $L^2
(\vOL,m)$ and the pure pointedness of the spectrum of $T$ follows. 

\smallskip

The last statement follows by a careful analysis of the steps in the proof of
$(ii)\Longrightarrow (i)$. 
\end{proof}

\begin{proof}[Proof of Theorem \ref{entropy}]
  For dynamical systems over $\ZZ$ it is well known that pure point spectrum
  implies vanishing of the metric entropy.   The reason is that pure point spectrum  implies that the system is measurably conjugated to a rotation on a compact Abelian group by the Halmos/von Neumann theorem. Such a rotation in turn has vanishing metric entropy. 
Also, for these systems a
  variational principle is well known relating topological and metric
  entropy.  Similar result can be shown for dynamical systems over $\RR^N$.
  (In the Euclidean case it is  an issue on how to define the entropy in the
  first place as we do not have a first return map.)  In fact, a variational
  principle can be found in the work \cite{TZ} of Tagi-Zade.

Given this the proof of the theorem proceeds along the following steps (see \cite{BLR} for details): 

Step 1: As $\gammahat$ is pure point, we have pure point dynamical spectrum by
Theorem \ref{characterization} and hence vanishing measurable entropy. 

Step 2: As $\Oo$ is uniquely ergodic, we obtain vanishing of the topological
entropy from Step 1 and the variational principle. 

Step 3: The topological entropy can be shown to be equal to the patch counting
entropy, which is the limit appearing in the theorem. 
\end{proof}

\begin{proof}[Proof of Theorem \ref{bombieritaylor}]
The proof is given in two steps. In the first step it is shown that uniform
convergence for the $c_n^\xi$ follows, once it is known to hold for certain
averages in a topological  Wiener/Wintner type  ergodic theorem.  In the second step, uniform convergence in
this Wiener/Wintner type theorem is then shown provided the eigenfunctions are
continuous.  For details  concerning this proof we refer to \cite{Le}. For a  general topological Wiener/Wintner ergodic theorem  and further references we refer to Section \ref{WW}.   
\end{proof}

\section{A word on symmetry}\label{Symmetry}
In this section we discuss a result on symmetries of $\gammahat$. The result
is a consequence of unique ergodicity and the closed formula. While it is
essentially a simple observation, we are not aware
of a reference. 

\begin{theorem} Let $\vL$ be Delone with (UPF) and (FLC).  Let $
  V:\RR^N\longrightarrow \RR^N$ be linear and isometric. If $\vOL$ is
  invariant under $V$, which means that  $V \vG =\{ Vx : x\in \vG\}\in \vOL$ for all
  $\vG\in \vOL$, then $\gammahat$ is invariant under $V$ as well i.e. 
$$ \int \varphi (V \xi) d\gammahat (\xi) = \int \varphi
(\xi)d\gammahat(\xi).$$
\end{theorem}
\begin{proof} As $V$ is linear, it is not hard to see that  the map 
$$ f\mapsto \int f(V\vG) d m(\vG)$$
is a translation invariant probability measure on $\vOL$. By unique
ergodicity, we then have
$$ \int f (V \vG) d(\vG) = \int f (\vG) d m (\vG)$$
for all $f\in C (\vOL)$. Define $\widetilde{V} \varphi (z) = \varphi (V z)$
for  $\varphi : \RR^N\longrightarrow \CC$. The closed formula for $\gammahat$
then shows that 
\begin{eqnarray*}
\gamma (\varphi) &=& \int \sum_{x,y\in \vG} \sigma (x) \varphi (x-y) dm
(\vG)\\
&=& \int \sum_{x,y\in V \vG} \sigma (x) \varphi (x-y) dm
(\vG)\\
&=& \int \sum_{x,y\in \vG} \sigma (V x) \varphi ( V (x-y)) dm
(\vG)\\
&=& \gamma (\widetilde{V} \varphi).
\end{eqnarray*}
Here, we used in the last step that $\int \sigma (V t ) dt =1$ and hence
$\gamma$ can be calculated with $\sigma \circ V$ as well as with $\sigma$. 
A short and direct calculation shows furthermore that
$$ (\widetilde{V} \widehat{\varphi}) (x) = \widehat{ \widetilde{V} \varphi}
(x).$$
Putting this together we obtain for all $\varphi \in \cS$
$$\gammahat (\widetilde{V} \varphi) =  \gamma (\widehat{ \widetilde{V}
  \varphi}) = \gamma ( \widetilde{V} \widehat{\varphi}  ) = \gamma (
  \widehat{\varphi}) = \gammahat(\varphi).$$
By density considerations, this formula then holds for all  $\varphi \in C_c  (\RR^N)$ and we obtain the statement. 
\end{proof}

\begin{remark} The proof does not use (FLC). It only uses unique ergodicity and  the closed formula. Accordingly, the result remains correct for uniquely ergodic situations without (FLC). 
\end{remark}

\section{A class of examples: Meyer sets}\label{Meyer}
There are two prominently studied classes of sets in the field of aperiodic order. Theses are 
 sets associated to primitive substitutions and Meyer sets. In this section we have a closer look at a special class of Meyer sets known as model sets.  Meyer sets can be thought of as very natural generalizations of lattices.  In fact, there are several characterizations of Meyer sets giving a precise meaning to this. 
Here, we shortly discuss an algebraic characterization due to Lagarias and then focus on a way to  create Meyer sets.  For further discussion and details we refer to \cite{Moo1,Moo2,Schl}

\begin{definition} A Delone set  $\vL$ in  $\RR^N$ is called Meyer set if  $\vL - \vL$ is uniformly discrete (and hence Delone as well). 
\end{definition}

As shown by Lagarias \cite{Lag} Meyer sets  in $\RR^N$ can be characterized by the following lattice like behavior (see \cite{BLM} for  a recent extension of Lagarias argument to certain locally compact Abelian groups as well).

\begin{theorem}[\cite{Lag}] A Delone set $\vL$ is Meyer if and only if there exists a finite set $F$  with
$$ \vL -\vL \subset \vL + F.$$
\end{theorem}

We now discuss how Meyer sets arise as projections from a higher dimensional lattice structure via so called cut and project schemes.

A cut and project scheme over $\RR^N$ consists of a locally compact Abelian group $H$,
called the internal space, and a lattice $\widetilde{L}$ in
$\RR^N\times H$ such that the canonical projection $\pi : \RR^N\times H
\longrightarrow \RR^N$ is one-to-one between $\tilde{L}$ and
$L:=\pi(\widetilde{L})$ and the image $\pi_{\rm
int}(\widetilde{L})$ of the canonical projection $\pi_{\rm int} :
\RR^N\times H\longrightarrow H$ is dense. Given these properties of
the projections $\pi$ and $\pi^{}_{\rm int}$, one can define the
$\star$-map $(.)^\star\!: L \longrightarrow H$ via $x^\star :=
\big( \pi^{}_{\rm int} \circ (\pi|_L)^{-1}\big) (x)$, where
$(\pi|_L)^{-1} (x) = \pi^{-1}(x)\cap\tilde{L}$, for all $x\in L$.

We summarize the features of a cut- and project scheme in the
following diagram:
\begin{equation*} \label{candp}
\begin{array}{cccccl}
    \RR^N & \xleftarrow{\,\;\;\pi\;\;\,} & \RR^N\times H &
        \xrightarrow{\;\pi^{}_{\rm int}\;} & H & \\
   \cup & & \cup & & \cup & \hspace*{-2ex} \mbox{\small dense} \\
    L & \xleftarrow{\; 1-1 \;} & \tilde{L} &
        \xrightarrow{\,\;\quad\;\,} & L^\star & \\
   {\scriptstyle \parallel} & & & & {\scriptstyle \parallel} \\
    L & & \hspace*{-38pt}
    \xrightarrow{\hspace*{47pt}\star\hspace*{47pt}}
    \hspace*{-38pt}& & L^\star
\end{array}
\end{equation*}
We will assume that the Haar measures on $\RR^N$ and on $H$ are chosen in
such a way that a fundamental domain of $\tilde{L}$ has measure $1$. Any cut and project scheme comes with a natural dynamical system $(\TT,\alpha')$. Here, $\TT := (\RR^N\times H) /\widetilde{L}$ and
$$ \alpha' : \RR^N \times \TT\longrightarrow \TT, \:\; \alpha_t'([s,h]) := [t +s,h].$$
By density of $L^\star$ this system is minimal. It can then easily be seen to be uniquely ergodic as well. 

\smallskip

Given a cut and project
scheme, we can associate to any $W\subset H$, called the window, the
set
\begin{equation*}
   \oplam(W) \; := \; \{ x\in L : x^\star \in W \}
\end{equation*}

The following two properties  of $\oplam (W)$ are well known.  We therefore only sketch the proof. 

\begin{prop} Let $(\RR^N, H, \widetilde{L})$ be a cut and project scheme. Let $W\subset H$ be given. 

(a) If the closure $\overline{W}$ of $W$ is compact, then $\oplam (W)$ is uniformly discrete. 

(b) If the interior $W^\circ$ of $W$ is not empty then $\oplam (W)$ is relatively dense. 
\end{prop}
\begin{proof} (a) Assume that there are points $x_n, y_n$ in $\oplam (W)$ with $x_n\neq y_n$ and  $x_n - y_n$ converging to $0$ for $n\to \infty$. These points come from points $(x_n, x_n^\star)$, $(y_n, y_n^\star)$ of the lattice. By assumption, $x_n^\star, y_n^\star \in W$. As $W$ is relatively compact, we can assume without loss of generality that  $(x_n^\star)$ and $(y_n^\star)$ are converging sequences with limits $\widetilde{x}$ and $\widetilde{y}$ respectively.  Consider now the sequence $ Z_n :=(x_n - y_n, x_n^\star - y_n^\star)$.   Our considerations show that the points $Z_n$ converge to  $Z=(0, \widetilde{x} - \widetilde{y})$.  Moreover, the points 
$Z_n$ belong to $\widetilde{L}$ as $\widetilde{L}$ is a lattice. Thus, $Z$ must belong to $\widetilde{L}$ as well. By the requirements on a cut and project scheme we infer  that $0 =\widetilde{x} - \widetilde{y}$. Hence, $(Z_n)$  is a sequence in the lattice converging to the origin. This is only possible, if $Z_n = (0,0)$ for large $n$. This contradicts $x_n\neq y_n$. 

(b)  Let $U$ be the open interior of $W$.  By definition $\RR^N \times H/ \widetilde{L}$ is compact. We can therefore find $S>0$ and $h_1,\ldots, h_n\in H$, such that 
$$ F_t := \cup (t + B_S)\times (h_i +U)$$
contains representatives of all elements in $\TT= (\RR^N \times H)/ \widetilde{L}$ for any $t\in \RR^N$. 
By density of $L^\star$ in $H$, we can assume without loss of generality that each $h_j$ belongs to $L^\star$ i.e. has the form $h_j = x_j^\star$ for some $x_j$ in $\RR^N$. Then
$$ ( t+ \cup (-x_j +B_S)) \times U$$
contains a representative of any element in $\TT$ for any $t\in \RR^N$.

If we now choose $R>0$ such that $B_R \supset \cup (-x_j +B_S)$ then any translate of $B_R$ will contain a point of $\oplam (U)$. 
\end{proof}

The proposition has the following consequence. 

\begin{coro} Let $(\RR^N, H, \widetilde{L})$ be a cut and project scheme and  $W\subset H$ relatively compact with non empty interior. Then, $\oplam (W)$ is Meyer. 
\end{coro} 
\begin{proof} By the previous proposition $\oplam (W)$ is Delone. As $\widetilde{L}$ is a lattice we have
$$\oplam (W) - \oplam (W) \subset \oplam ( W - W).$$
As $W$ is relatively compact, so is $W- W$ and we infer  from (a) of the previous proposition that $\oplam (W) - \oplam (W) $ is uniformly discrete. 
\end{proof}

A set of the form $t+ \oplam(W)$ is called model set if the  window $W$ is relatively compact with
nonempty interior.  The following remarkable converse of the previous corollary  holds \cite{Mey,Moo1}. 

\begin{theorem} A subset $\vL$ of $\RR^N$ is Meyer if and only if it is a subset of a model set. 
\end{theorem}

A model set is called {\em regular\/} if $\partial W$ has Haar measure $0$ in $H$.
For $\vL=\oplam (W)$ with relative compact $W$ which is the closure of its interior, there is s strong connection between the dynamical system $\Oo$ and the canonical dynamical system 
$(\TT,\alpha')$ introduced above.  This connection  is given as follows (see
Proposition 7 in \cite{BLM} for the statement given next and
\cite{Schl}, \cite{BHP} for earlier versions of the same type of
result).

\begin{prop}\label{blm-prop} 
There exists a continuous $\RR^N$-map $\beta : \varOmega(\vL)
\longrightarrow \TT$ with the property that $\beta(\Gamma) = (t,h) +
\widetilde{L}$ if and only if $t + \oplam (W^\circ - h) \subset
\Gamma \subset t+ \oplam(W - h)$.
\end{prop}

Using this proposition (or similar results) it is possible to conclude properties of $\Oo$ from properties of $(\TT,\alpha)$.  If $\varLambda$ is regular, then  the map $\beta$ is almost everywhere $1:1$ by the previous proposition. Thus,  we can easily infer the following properties of $\Oo$ from the corresponding properties of  $(\TT,\alpha')$  (see e.g.  \cite{Schl,BLM}):
\begin{itemize}
\item $\Oo$ is uniquely ergodic.

\item $\Oo$ has pure point dynamical spectrum.

\item All eigenfunctions of $\Oo$ are continuous. 
\end{itemize}
In particular, we obtain pure point diffraction \cite{Hof,Schl}.  In this case, one can  calculate explicitely the diffraction measure $\gammahat$ \cite{Hof,Schl}.
For $x\in \RR^N$ and $k\in \widehat{\RR^N} = \RR^n$, we set
$$ (k,x):= e^{i k x}.$$
We need the  dual lattice $\Lperp$ of $\widetilde{L}$ given by
$$\Lperp :=\{ (k,u) \in\widehat{\RR^N}\times \widehat{H} : k(l) u(l^\star) =1\;\:\;\mbox{for all
$(l,l^\star) \in \widetilde{L}$}\}.$$ Let $\Lnull$ be the set of all
$k\in \widehat{\RR^N}$ for which there exists $u\in \widehat{H}$ with
$(k,u)\in \Lperp$. As $\pi_2 (\widetilde{L})$ is dense in $H$, we
easily infer that $(k,u), (k,u')\in \Lperp$ implies $u=u'$. Thus,
there exists a unique map $\star : \Lnull\longrightarrow \widehat{H}$
such that
$$\tau: \Lnull \longrightarrow \Lperp, \,\:\; k\mapsto (k,k^\star)$$
is bijective.  Then, the diffraction measure $\gammahat$ is given by
$$  \gammahat = \sum_{k\in \Lnull} A_k \delta_k,$$
where $A_k = |\int_W (k^\star,y)dy |^2$. We will shortly sketch a proof  based on Theorem \ref{bombieritaylor} above, see  Lenz/Strungaru  \cite{LS} as well. 
For $n\in \NN$ and $k\in \Lnull$ consider the function
$$ c^k_n : \vOL\longrightarrow \CC, \;\;\mbox{by}\;\: c^k_n
({\Gamma} ) :=\frac{1}{\left|B_n\right|} \sum_{x\in
{\Gamma}\cap B_n} \overline{(k,x)}.$$ By  Theorem  \ref{bombieritaylor} above and the stated properties of $\Oo$, we know that the $c^k_n$ converge  uniformly to a function $c^k$ and $|c^k|^2$
equals the coefficient $A_k$. Thus, it remains to
calculate the limit of the function $c^k_n$. 

Now, by Proposition \ref{blm-prop},
$\Gamma$ has the form $\Gamma = t + \widetilde{\Gamma}$ with
$\oplam (-h +W^\circ) \subset
\widetilde{\Gamma} \subset \oplam (-h + W)$ for suitable $t\in \RR^N$ and
$h\in H$. For
$k\in \Lnull$ and $x\in L$ we have by definition $$(k,x) =
\overline{(k^\star,x^\star)}.$$
 Moreover, denoting the character $(k,k^\star)\in \widehat{\RR^N} \times \widehat{H}$ by $\tau (k)$ we find by direct calculation
 $(\tau(k), \beta(\Gamma)) = (k,t) (k^\star,h )$ and hence
$$ \overline{(k,t)} = (k^\star,h ) \overline{  (\tau(k), \beta(\Gamma))}.$$
Combining all of this we obtain
$$ \overline{(k,t + x)} =\overline{ (\tau(k),\beta(\Gamma)) }.$$
Thus, the term of interest is given by
$$\frac{ \overline{ (\tau(k),\beta(\Gamma))} }{|B_n|} \sum_{x
\in\widetilde{\Gamma}\cap (-t + B_n)} (k^\star,x^\star + h) .$$

By uniform distribution
\cite{Moo3} this converges  to
$$ c^k (\Gamma) = \overline{ (\tau(k),\beta(\Gamma))   } \int_{-h + W}
(k^\star,y + h )  dy = \overline{ (k,k^\star) (\beta(\Gamma))} \int_{W}
(k^\star,y  ) dy     .$$

Thus, $|c^k (\Gamma)|^2$ is equal to $$A_k := \left|\int_W (k^\star,y)  dy\right|^2.$$

\bigskip

\medskip

\bigskip

\section{Uniform Wiener/Wintner type theorems}\label{WW}
In this section we present a Wiener/Wintner type theorem for actions of $\RR^N$.  A theorem of this kind lies at the heart of our proof of theorem \ref{bombieritaylor} sketched above. The  theorem given here is new. It generalizes a main result of \cite{Le},   which in turn generalizes a result Robinson \cite{Rob}. At the same time our result extends a result of Walters from actions of $\NN$ to actions of  $\RR^N$.   Our proof is essentially an extension of ideas developed in \cite{Le}. For this reason we only sketch it.  The result is valid for general topological dynamical systems and not only systems coming from Delone sets. Accordingly, we work in a slightly more general setting here than in the rest of the paper.  For related results we also refer to work of Assani \cite{Ass}.

\medskip

As usual $(\varOmega,\alpha)$ is called a topological dynamical system  over $\RR^N$ if $\varOmega$ is a compact topological space and $\alpha : \RR^N \times \varOmega \rightarrow \varOmega$ is a continuous action of $\RR^N$ on $\varOmega$. Denote the set of continuous functions on $\varOmega$ by $C(\varOmega)$.   Let $S$ denote the  unit circle in the complex plane. Given a topological dynamical system $(\varOmega, \alpha)$ a continuous map $\varphi : \RR^N \times \varOmega \longrightarrow S$ is called a cocycle if
$$ \varphi (x + y, \omega) = \varphi (x,\alpha_y (\omega)) \varphi (y,\omega)$$
for all $x,y\in \RR^N$ and $\omega \in \varOmega$.   Any character $(\xi, \cdot) : \RR^N \longrightarrow S, (\xi,x):= e^{i\xi x}$ (where $\xi\in \RR^N$) induces a cocycle viz $\varphi_\xi (x,\omega) := (\xi,x)$. 

 Let now a dynamical system  $(\varOmega, \alpha)$ together with an ergodic probability  measure $m$ be given. Each cocycle  on this system gives then  rise to a unitary representation $\Uphi$ of $\RR^N$ on $L^2 (\varOmega,m)$ via
$$ \Uphi_t f (\omega) := \varphi (t,\omega) f(\alpha_{-t} \omega)$$
for $t\in \RR^N$. By ergodicity and the usual arguments, the  subspace of  solutions to $ \Uphi_t f  = f$ for all $t\in \RR^N$ is one dimensional or trivial.  
Let $\Pphi$ be the projection onto this subspace.

Define for $n\in \NN$ the map
$$ \ani : C(\varOmega)\longrightarrow C (\varOmega), \ani (f) (\omega):= \frac{1}{|B_n|} \int_{B_n} (\Uphi_t f)(\omega) dt.$$

\begin{theorem} Let a topological dynamical system $(\varOmega, \alpha)$ and a continuous  cocycle $\varphi$ be given.  Let $f\in C(\varOmega)$ be given. The following assertions are equivalent. 

\begin{itemize}

\item{(i)} $\Pphi f$ is continuous (i.e. there exists $g\in C(\varOmega)$ with $ g = \Pphi f$ in $L^2$ sense and  $\varphi (t,\omega)  g(\alpha_{-t} \omega)  = g(\omega)$ for all $t\in \RR^N$ and $\omega \in \varOmega$). 

\item{(ii)} The sequence $(\ani (f))$ converges uniformly (i.e. with respect to the supremum norm) to $\Pphi f$. 
\end{itemize}
\end{theorem}
\begin{proof}   As mentioned already this can be shown using the method developed in 
\cite{Le}. We sketch the proof.

By von Neumann ergodic theorem, the sequence $(\ani (f))$ converges to  $\Pphi f$ in the $L^2$ sense. 

(ii)$\Longrightarrow$ (i): This is clear as each $\ani (f)$ is continuous. 

(i) $\Longrightarrow$ (ii): Let $\varepsilon >0$ be given. By $L^2$-convergence of $(\ani (f))$ to $\Pphi f$, there exists an $N\in \NN$ such that the measure of  
$$\varOmega_N := \{ \omega \in \varOmega : | \ani (f) (\omega) - \Pphi f (\omega)| \geq \varepsilon \}$$
is smaller than $\varepsilon$.  By (i) and continuity of the $\ani(f)$, the set $\varOmega_N$ is closed and hence compact. This will be crucial. 

For $n$ large compared to $N$, $\ani (f)$ and $\aNi ( \ani(f))$ become  arbitrarily  close to each other. It therefore suffices to consider $\aNi (\ani(f))$.  By Fubini's theorem this is equal to  $\ani (\aNi (f)) $. Let $\chi_N$ be the characteristic function of $\varOmega_N$.    Then, 
\begin{equation}\label{stern} \ani (\aNi (f))  = \ani ( (1-\chi_N) \aNi (f)) + \ani (\chi_N \aNi (f)).
\end{equation}
By unique ergodicity and compactness of $\varOmega_N$, we have 
$$\liminf_{n\to \infty} \frac{1}{|B_n|} \int_{B_n}  \chi_N (\alpha_{-t} \omega) \leq m( \varOmega_N) \leq \varepsilon$$
uniformly in $\omega \in \varOmega$.  This makes the second term  in \ref{stern} small in the supremum norm.  On the other hand,  in the first term  in \ref{stern} we can replace  $\aNi (f)$ by $\Pphi (f)$ and this will be a small error by the very definition of $\varOmega_N$.  These considerations show that  $\ani (\aNi (f))$ is close to $\ani (\Pphi (f))$ for large $n$. The latter, however, equals $\Pphi (f)$ by definition of $\Pphi$. This finishes the proof. 
\end{proof}

\begin{remark} (a) Note that  (ii) in the above theorem contains trivially the case that   $\Pphi f=0$.  Thus,  (ii) covers both the situation that there does not exist an eigenfunction of $\Uphi$ to the eigenvalue $1$ and the situation that there exists a continuous eigenfunction. Theses two cases are investigated separately by Walters \cite{Wal2} for actions of $\NN$. For cocycles coming from characters these  cases are investigated for  actions of $\NN$ and $\RR^N$  by Robinson  \cite{Rob}.  For cocycles coming from characters and actions of locally compact Abelian group the result above is given in \cite{Le}.  

(b) The result above is stated and proved for $\RR^N$. The crucial ingredient, however, is the validity of a von Neumann Ergodic Theorem. Such a theorem is known for locally compact Abelian groups and for various semigroups e.g. $\NN$ (see e.g. \cite{Krengel}). Thus, the proofs and results carry over to these situations as well. 
\end{remark}

\section{Further remarks and open questions}\label{Further}
We have discussed a framework for diffraction based on Delone sets with (FLC) in Euclidean space.   This is a natural framework when  one wants to preserve the connection  between diffraction and geometry/combinatorics.    One may ask, however, how necessary these assumptions really are. This is not only of abstract mathematical interest in helping understanding the assumptions. It is also relevant from the point of view of modeling physical substances.  For this purpose one may well argue that more general point sets should be admissible or, even, that point sets are too restrictive altogether. Accordingly, various generalizations have been considered. 

It turns out that $\RR^N$ can be replaced by an arbitrary locally compact Abelian group when dealing with diffraction for (FLC) sets and model sets.  This is carried out in work of Schlottmann \cite{Schl}. Likewise one may consider more general point sets as discussed (in Euclidean space) by Gou\'{e}r\'{e} \cite{Gou1,Gou2}. In fact, one can leave the framework of point sets altogether and work with measures instead. This is studied (on locally compact Abelian groups) by Baake/Lenz \cite{BL,BL2}, Lenz/Richard \cite{LR} and Lenz/Strungaru \cite{LS}. 

Even within the framework studied above various questions and issues present themselves. Here, we would like to mention the following questions (see the survey article of Lagarias \cite{Lag2} as well). 

\medskip

The discussion above gives the following chain of inclusions:

Lattices $\subset$ Regular model sets $\subset$  Meyer sets with pure point diffraction $\subset $ Delone sets with (FLC) and pure point diffraction $\subset$ Delone sets with (FLC) and a relatively dense set of Bragg peaks. 

\smallskip

\textbf{Question.} How far are these inclusions from being strict or put differently, how can one  characterize each of these classes of sets within the next bigger class?

\smallskip

A natural issue in this context is the following. 

\textbf{Question.} Does existence of (pure) point diffraction together with some further conditions like (FLC) and repetitivity already imply the Meyer property? 

\smallskip

It seems that the only results in this direction are  proven within the context of primitive substitutions. Lee/Solomyak   \cite{LSo} show that the Meyer property follows for primitive substitutions with pure point spectrum.  Lee \cite{Lee} then shows that for primitive substitutions pure point diffraction is in fact equivalent to being a  model set. 

On the other hand, by recent results of Strungaru \cite{Str}, the Meyer property already implies existence of a relative dense set of Bragg peaks. 

\smallskip

\textbf{Question.}  What is the significance of  a relatively dense set of Bragg peaks?

\medskip

Finally, we note that our discussion of model set heavily relied on the assumption of regularity i.e. vanishing measure of the boundary of the window. 

\textbf{Question.} What can one say about model sets with a thick boundary?

\subsection*{Acknowledgments}  The author would like to thank Idris Assani for the invitation to the  ergodic theory workshop in spring  2007 in Chapel Hill and Karl Petersen for  discussions at this workshop. These lead to the  material presented in Section \ref{WW}. The author  would also like to take the opportunity to thank Michael Baake, Robert V. Moody, Christoph Richard and Nicolae Strungaru for a both  enjoyable and educational collaboration on the topic of  aperiodic order.  Partial support from the German Science Foundation (DFG) is gratefully acknowledged.

\bigskip
\bigskip

\bibliographystyle{amsalpha}

\end{document}